\documentclass{article}

\usepackage{amssymb}
\usepackage{amsthm}

\newtheorem{thm}{Theorem}
\newtheorem{prop}{Proposition}

\begin{document}

\title{Spectra of symmetric powers of graphs and the Weisfeiler-Lehman refinements}
\author{Afredo Alzaga  \\Departamento de Matem\'{a}tica \\ Universidad Nacional del Sur (UNS) \\
        Bah\'{i}a Blanca, Argentina
 \and Rodrigo Iglesias  \\Departamento de Matem\'{a}tica \\ Universidad Nacional del Sur (UNS) \\
        Bah\'{i}a Blanca, Argentina
 \and  Ricardo Pignol  \\Departamento de Matem\'{a}tica \\ Universidad Nacional del Sur (UNS) \\
        Bah\'{i}a Blanca, Argentina}  
\maketitle

\begin{abstract}
The $k$-th power  of a $n$-vertex graph $X$ is the iterated cartesian product of $X$ with itself. The $k$-th symmetric power of $X$ is the quotient graph of certain subgraph of its $k$-th power by the natural action of the symmetric group. It is natural to ask if the spectrum of the $k$-th power --or the spectrum of the $k$-th symmetric power-- is a complete graph invariant for small values of $k$, for example, for $k=O(1)$ or  $k=O(\log n)$.

In this paper, we answer this question in the negative: we prove that if the well known $2k$-dimensional Weisfeiler-Lehman method fails to distinguish two given graphs, then  their $k$-th powers --and their $k$-th symmetric powers-- are cospectral.  As it is well known, there are pairs of  non-isomorphic $n$-vertex graphs which are not
distinguished by the $k$-dim WL method, even for $k=\Omega(n)$. In particular, this shows that for each $k$, there are pairs of non-isomorphic $n$-vertex graphs with cospectral $k$-th (symmetric) powers.
\end{abstract}

\section{Introduction}

Many fundamental graph invariants arise from the study of random walks of a particle on a graph.
 Most of these invariants
 can be described in terms of the spectrum of the adjacency or the Laplacian matrix.
Since the graph spectrum fails to distinguish many  non-isomorphic
graphs, it is interesting to study the properties of walks
 (or quantum walks)  of $k$ particles, as a means to construct more powerful invariants.

This led Audenaert \emph{et al} \cite{Godsil} to define the $k$-th
symmetric power $X^{\{k\}}$ of a graph $X$: each vertex of
$X^{\{k\}}$ represents a $k$-subset of vertices of $X$, and two
$k$-subsets are joined if and only if their symmetric difference
is an edge of $X$. They show that the spectra of these graphs  is
a family of invariants stronger than the ordinary graph spectra.
For $k=2$, they provide examples of cospectral graphs $X$ and $Y$
such that $X^{\{2\}}$ and $Y^{\{2\}}$ are not cospectral. On the
other hand, they prove that if $X$ and $Y$ are strongly-regular
cospectral graphs then $X^{\{2\}}$ and $Y^{\{2\}}$ are cospectral.
For $k=3$, the authors reported computational evidence suggesting
that the spectra of the symmetric cube may be a strong invariant.
They did not find any pair of non-isomorphic graphs with
cospectral $3$-symmetric powers, upon inspection of all strongly
regular graphs of up to $36$ vertices.

In this paper we prove that for each $k$ there are pairs of non-isomorphic graphs such that their $k$-th symmetric
 powers are cospectral
by showing how these invariants are related to the well
known $k$-dimensional Weisfeiler-Lehman (WL) algorithm.

The automorphism group of the graph acts on the set of $k$-tuples
of vertices. The $k$-WL method is a combinatorial algorithm  that
attempts to find the associated orbit partition (see for example \cite{Cai}, \cite{Grohe 2000}). It starts by
classifying the $k$-tuples according to the isomorphism type of
their induced graphs, then  an iteration is performed attaching to
the previous color of a $k$-tuple, the multiset of colors of the
the neighboring $k$-tuples. In this way, the partition of the
$k$-tuples is refined in each step until a stable partition is
reached. The multiset of colors of the stable partition is a graph
invariant.

Our main result is the following theorem.
\begin{thm} \label{Main theorem}
If the $2k$-dim Weisfeiler-Lehman algorithm fails to distinguish two given graphs, then their
$k$-th symmetric powers are cospectral.
\end{thm}
In fact, the result remains true if we consider $k$-th powers of
graphs (associated to  walks of $k$ labelled particles), instead
of symmetric powers.

In \cite{Cai},  Cai, Immerman and F\"urer showed how to construct
pairs of non-isomorphic $n$-vertex graphs which are not
distinguished by the $k$-WL method, even for $k=\Omega(n)$. Then,
our result implies that
\begin{thm} \label{Main corollary}
If we require the $k$-th symmetric power spectrum to determine all
$n$-vertex graphs then, necessarily, $k=\Omega(n)$.
\end{thm}

Nevertheless, the spectrum of the $k$-th power of a graph  is a
strong invariant with remarkable computational features. Since it
is determined by the characteristic polynomial of a matrix of $0$s
and $1$s  of polynomial size (for fixed $k$), it   can be computed
in polylogarithmic time by a randomized parallel algorithm. This
contrasts with the inherently sequential nature of the $k$-dim WL
algorithm; in \cite{Grohe}, Grohe proved that finding the $k$-dim
WL stable partition is a $P$-complete problem ($k \geq 2$).
This suggest that the $2k$-WL method ($k \geq 1$) is strictly
more powerful than the $k$-th power (or $k$-th symmetric power)
spectrum, since the complexity class $RNC$ is expected to be strictly contained in $P$.

Besides power graph spectra, there are other families of  graph
invariants in the literature for which it is not known whether
they distinguish any pair of non-isomorphic graphs or not. As it
turns out,   the WL-refinements provide a natural benchmark to
compare other graph invariants and it is reasonable to expect that
arguments of the kind we use in this work would show the
limitations of some of them.

The paper is organized as follows. In Section \ref{Powers of
graphs} we define the $k$-th power $X^{k}$ and the $k$-th
symmetric power $X^{\{k\}}$ of a graph $X$. In Section
\ref{Quotient graphs} we recall the general notion of quotient of
a graph by the action of a group, and we  describe the  $k$-th
symmetric power as a quotient of the restricted $k$-th power
$X^{(k)}$. For later use, we prove some formulas  concerning the
walk generating function of quotient graphs. In Section \ref{The
Weisfeiler-Lehman algorithm} we define precisely the
$k$-Weisfeiler-Lehman algorithm. The heart of the proof of Theorem
$1$ is in Section \ref{Spectrum of the $k$-th power vs.  the
$2k$-WL coloring}. Essentially, we show that the $2k$-WL method is
stronger than the spectra of the $k$-th power $X^{k}$. Since the
idea of the proof is easier to exhibit in the  case $k=1$, we
write this special case separately in Section  \ref{Graph spectra
is weaker than the $2$-WL refinement}. Finally, the proof of
Theorem $1$ is given in  Section \ref{Proof of Theorem 1}, by
passing to the quotient $X^{\{k\}}$. In order to achieve this, we
 exploit the  structure of the set of $k$-tuples and the
formulas for quotient graphs presented in Section \ref{Quotient
graphs}.

\section{Powers of graphs}
\label{Powers of graphs}

In this section we present the notion of the  $k$-th symmetric
power of a graph, as introduced  in \cite{Godsil}, and some other related
constructions.

Throught the paper, a \textit{graph} $G$ is a finite set $V$ of vertices toghether
with a set $E$ of unordered pairs $(v,w)$ of vertices with $v \neq w$.
We denote by $A_G$ the \textit{adjacency matrix of G}. Since we do not assume an order on $V$, we consider $A_G$ as a function $A_{G}: V \times V \rightarrow \mathbb{Z}$,  defined by
 $A_{G}(v,w)=1$ if $(v,w) \in E$, and
$A_{G}(v,w)=0$ otherwise.

A \textit{k-tuple} $(i_1...i_k)$ of vertices is a function from $\{1,...,k\}$ to $V$. Let $\mathcal{U}_k$ be the set of  all $k$-tuples and let $\mathcal{D}_k \subset \mathcal{U}_k$ denote the
set of those $k$-tuples of pairwise distinct vertices.
The symmetric group $S_k$ acts naturally on $\mathcal{D}_k$ by
$
\sigma (i_1...i_k) = (i_{\sigma^{-1}(1)}...i_{\sigma^{-1}(k)}),
$
for $\sigma \in S_k$.  The orbits are identified with the \textit{k-subsets} of vertices.

The $k$-th symmetric power of $G$, denoted by $G^{\{k\}}$, has the
$k$-subsets of $V$ as its vertices; two $k$-subsets  are adjacent
if their symmetric difference --elements in their union but not in
their intersection-- is an edge of $G$. The picture behind this
construction  is borrowed from the physical realm: start with $k$
undistinguishable particles occupying $k$ different vertices of
$G$ and consider the dynamics of a walk through the graph in
which, for each step, any single particle is allowed to move to an
unoccupied adjacent vertex.  In this way, a $k$-walk on $G$
corresponds to a $1$-walk on $G^{\{k\}}$.  The connection between
symmetric powers and quantum mechanics exchange Hamiltonians is
further explored in [2].

Likewise, one can define the \emph{cartesian product}  $G \times
H$ of two graphs as follows:
$$
 A_{G\times H}(i_1 i_2 ,
j_1 j_2 ) = \left\{
\begin{array}{ll}
1 \ & \mbox{if } \ A_G(i_1,j_1) =1 \ \mbox{and }\   i_2 = j_2 \ \\
    & \mbox{or else } \ A_G(i_2,j_2) =1 \ \mbox{and }\   i_1 = j_1 \ \\
0 &\mbox{otherwise}
\end{array}\right.
$$
 The $k$-th \textit{power} $G^k$ of a graph is
 defined as the iterated cartesian product of $G$ with itself. The set of its vertices is
 $\mathcal{U}_k$ and its adjacency matrix $A_{G^k}$ is given by:
$$
 A_{G^k}(i_1 i_2 \ldots i_k ,
j_1 j_2 \ldots j_k) = \left\{
\begin{array}{ll}
1 &\mbox{if there exists }\  u \in \{1,\ldots,k\} \  \mbox{such
that
}\\
&A_G(i_u,j_u) =1 \ \mbox{and }\   i_l = j_l \  \mbox{for } \  l\neq u \\
0 &\mbox{otherwise}
\end{array}\right.
$$
In the physical cartoon of the particles, the $k$-th power
correspond to the  situation in which the $k$ particles are labeled,
and more than one particle is allowed to occupy the same vertex at the same time.

Given a graph $G$, the \textit{walk generating function} of $G$ is
the  power series
$$
 \sum_{r=0}^{\infty}{t^{r} (A_{G})^{r}}
$$
The coefficient of $t^r$ in the $(i,j)$-entry counts the number of paths of
 length $r$ from the vertex $i$ to the vertex $j$.
See \cite{Godsil} for further properties. The trace of the walk
generating function is a graph invariant, and we denote it by
$$
F(G,t)=Tr \sum_{r=0}^{\infty}{t^{r} (A_{G})^{r}}
$$
Since the spectrum of two matrices $A$ and $B$ coincides if and
only if $Tr(A^r)=Tr(B^r)$ for all $r$, two graphs $G$ and $H$ are
cospectral if and only if $F(G,t)=F(H,t)$. In particular, they
cannot be distinguished by the spectrum of their $k$-th symmetric
powers  if and only if $F(G^{\{k\}},t)=F(H^{\{k\}},t)$.

\section{Quotient graphs}
\label{Quotient graphs}

The $k$-th symmetric power $G^{\{k\}}$ can be constructed from $G^k$
in two steps. First, we cut $G^k$, deleting all those vertices
which are not in $\mathcal{D}_k$. In this way we obtain the
\textit{restricted $k$-th power}, denoted  by $G^{(k)}$, defined
as the subgraph of $G^{\{k\}}$ whose vertices are the $k$-tuples in
$\mathcal{D}_k$. Second, we take  the quotient of $G^{\{k\}}$ by
the natural action of $S_k$ on the restricted $k$-th power $G^{(k)}$.

Let us give the general definition of a quotient graph and discuss some properties.
Given a graph $X$ and a group
$\Gamma$ acting on $X$ by automorphisms, the quotient $X/ \Gamma$ is a directed graph, in general with multiple edges and loops, defined as follows. The vertices of $X/ \Gamma$ are the orbits of the vertices of $X$, and given two orbits $U$ and $W$, there are as many arrows from $U$ to $W$ as edges in $X$ connecting a fixed element $u \in U$ with vertices in $W$.

We are interested in the case where this quotient has no loops and
no multiple edges; we say that the quotient $X/\Gamma$ is
\textit{simply laced} if
\begin{enumerate}
    \item $(u,v) \in E$ implies that $u$ and $v$ are not in the same orbit.
    \item $(u,v) \in E$ and $(u,w) \in E$  implies that $v$ and $w$ are not in the same orbit.
\end{enumerate}
If $X/\Gamma$ is simply laced, we can consider it an ordinary graph, where $(U,W)$ is an edge if and only if there is an arrow in $X/\Gamma$ connecting them.

In the simply laced case, every path on $X/\Gamma$ can be lifted to an essencially unique  path on $X$. This fact simplifies the task of path-counting, and allows to derive  a simple formula for the walk generating function of a quotient graph. We apply it to the  symmetric power $G^{\{k\}}$ to obtain a formula that will be useful later.

\begin{prop} \label{path lifting}
Let $X$ be a graph, $X/\Gamma$ a simply laced quotient, and let $U$ and $W$ be two orbits. Then, the $r$-th power of the adjacency matrix of $X/\Gamma$ is given by
$$
A_{X/\Gamma}^{r}(U,W) = \frac{1}{\left|U\right|} \sum_{u \in U} \sum_{w \in W} A_{X}^{r}(u,w)
$$
\end{prop}

\textit{Proof:}
The entry $A_{X/\Gamma}^{r}(U,W)$ equals the number of paths of length $r$ on $X/\Gamma$
from $U$ to $W$. Fix an element $u_0 \in U$ and let $V_0,V_1,V_2,...,V_r$ be a path of
 length $r$ on $X/\Gamma$, with $U=V_0$ and $V_r=W$. Since there is at most one edge
 in $X$ connecting a vertex in $X$ to a vertex in a different orbit, there is a unique
  path $v_0,v_1,v_2,...,v_r$ in $G$ such that $v_0 = u_0$ and $v_j \in V_j$ for $0 \leq j\leq r$. Then,
$$
A_{X/\Gamma}^{r}(U,W) = \sum_{w \in W} A_{X}^{r}(u_0,w)
$$
The set of paths of length $r$ from $u_0$ to  $W$ is carried bijectively to the set of paths from any $u\in U$ to $W$ via some automorphism in $\Gamma$. Then, the sum $$\sum_{w \in W} A_{X}^{r}(u,w)$$ does not depend on $u$, and this proves the formula of the proposition.  \qed

Observe that this formula implies that  if $X/\Gamma$ is a
connected, simply laced quotient, then all the orbits have the
same size.

Let $M_{X/\Gamma}$ be the matrix  with rows and columns indexed by
the vertices of $X$, defined by
$$
M_{X/\Gamma}(v,w)=
\left\{
\begin{array}{ll}
\left|U\right| \ & \mbox{if } \ v \ \mbox{and }\   w \ \mbox{are in the same orbit} \ U  \\

0 &\mbox{otherwise}
\end{array}
\right.
$$
From Prop. \ref{path lifting} it follows:

\begin{prop} \label{Trace formula}
Let  $X/\Gamma$ be  simply laced quotient, and let $M_{X/\Gamma}$ be defined as above. Then,
$$
Tr (A_{X/\Gamma}^{r}) =  Tr (A_{X}^{r}M_{X/\Gamma}).
$$
\end{prop}

Now we set $X = G^{(k)}$ and $\Gamma = S_k$, acting in the natural
way on $G^{(k)}$. The  quotient $G^{(k)}/ S_k$ is isomorphic to
the $k$-th symmetric power $G^{\{k\}}$, and it is easily seen to
be a simply laced quotient. In this case, the matrix
$M_{X/\Gamma}$ is the matrix $M_{k}$, with rows and columns
indexed  by $k$-tuples in $\mathcal{D}_k$,  given by
$$
M_{k}(i_1...i_k,j_1...j_k)=
\left\{
\begin{array}{ll}
k! \ & \mbox{if } \ \{i_1...i_k\} \ \mbox{and }\   \{j_1...j_k\} \ \mbox{are equal as sets}   \\
    0 &\mbox{otherwise}
\end{array}
\right.
$$
From Prop. \ref{Trace formula} we obtain:
\begin{prop} \label{Quotient}
Let $G^{(k)}$ and  $G^{\{k\}}$  be the restricted $k$-th power and the $k$-th symmetric power of a graph $G$, respectively. Let $M_k$ be the matrix defined as above. Then,
$$
Tr (A_{G^{\{k\}}}^{r}) = Tr (A_{G^{(k)}}^{r}M_k)
$$
\end{prop}

\section{The Weisfeiler-Lehman algorithm}
\label{The Weisfeiler-Lehman algorithm}

A natural approach to graph isomorphism testing is to develop
algorithms to compute the vertex orbits of the automorphism group
of a  graph. In particular, if  the orbits of the union of two
graphs are known, one can decide if there is an isomorphism
between them. As a first approximation to the orbit partition of a
given graph, one can assign different colors to the vertices
according to their degrees. We can refine this partition
iteratively, by attaching to the previous color of a  vertex, the
multiset of colors of its neighbors. After at most $n=\left| V
\right|$ steps, the partition stabilizes. For most graphs, this
method distinguishes all the vertices \cite{Babai}, but it does
not work in general. For example, it clearly fails if the vertex
degrees are all equal to each other.

A more powerful method,  generalizing the previous one, is
obtained by coloring the $k$-tuples of vertices (single vertices
are implicit as $k$ repetitions of the same vertex). We start
classifying the $k$-tuples according to the isomorphism type of
their induced labelled graphs. Next, we apply an iteration
attaching to the previous color of a $k$-tuple, the multiset of
colors of the the neighboring $k$-tuples. This is  the so called
$k$-dimensional Weisfeiler-Lehman refinement. For fixed $k \geq 1$
the partition of the  $k$-tuples is no longer refined after $n^k$
steps, so the algorithm runs in polynomial time.

This type of combinatorial methods have been investigated since
the seventies, and for some time there was  hope in solving the
graph isomorphism problem provided that  $k=O(\log n)$ or
$k=O(1)$. In \cite{Cai}, Cai, Immermann and Fuhrer, disposed of
such conjectures; they proved that, for large $n$, $k$ must be
greater than $c n$ for some constant $c$, if we require the $k$-WL
refinement to reach the orbit partition of any $n$-vertex graph.
Despite of this limitation, the method works with $k$ constant
when restricted to some important
 families, such as planar or bounded genus graphs \cite{Grohe 2000}.

Let us define the $k$-WL algorithm more precisely.
Let $G$ be a graph and $V$ its set of
vertices. We define an equivalence relation on the set of $k$-tuples:
we say that $(i_1 \dots i_k)$ and $(j_1 \dots j_k)$ are equivalent if
\begin{enumerate}
    \item $i_{l}=i_{l'}$ if and only if $j_{l}=j_{l'}$
    \item $(i_{l},i_{l'}) \in E$ if and only if $(j_{l},j_{l'}) \in E$
\end{enumerate}

 We define the  \emph{type} $\:tp\:(i_1 \dots i_k)$ of a
$k$-tuple as its equivalence class.
Let $S_1$ be the set of all different types of $k$-tuples.
 This is the  initial set of colors. We define the
set $S$ of colors by
$$
S = \bigcup_{k=1}^{\infty}S_k
$$
where elements of $S_{r+1}$ are finite sequences or finite
multisets of elements of $ \bigcup_{k=0}^{r}S_k $. In practice, it
suffices to work with as many colors as $k$-tuples: in order to
preserve the length of their names, the colors can be relabelled
in each round (using  a rule not depending on $G$). Nevertheless,
this relabelling plays no role in our arguments.

We denote   the color assignment of the $k$-WL iteration in its $r$-th round,
applied to the graph $G$,
by $W_{G,k}^r : \mathcal{U}_k \rightarrow S$. Evaluated at the $k$-tuple $(i_1\dots i_k)$ it
gives the color $W_{G,k}^{r}(i_1\dots i_k) \in S$.
Initially, for $r=1$, it is defined by
$$W_{G,k}^{1}(i_1\dots i_k) =
\:tp\:(i_1\dots i_k).$$
The iteration is given by
\begin{equation} \label{Definition of k-WL}
W_{G,k}^{r+1}(i_1\dots i_k) = \sum_{m \in V} \left(\:tp\:(i_1\dots i_k \: m),
S_{G,k}^r(i_1 \dots i_k \: m)\right)
\end{equation}
where $S_{G,k}^r(i_1 \dots i_k \: m)$ is the sequence
$$
 \left(W_{G,k}^{r}(i_1 \dots \:
m),\dots, W_{G,k}^{r}(i_1\dots m \dots i_k),\dots,
W_{G,k}^{r}(m \dots i_k)\right).
$$
 The  summation symbol in (\ref{Definition of k-WL}) must be interpreted as  a formal sum,
 so that it denotes
 a multiset. For example, if $x_1 = x_3 = x_4 = a$ and $x_2 = x_5  = b$,

then $\sum_{i=1}^{5}x_i$ is the multiset
 $\{ a,a,a,b,b   \}$.

 For each round,
a certain number of different colors is attained. We say that the coloring
scheme \textit{stabilizes} in the $r$-th round if the number of different
colors does not increase in the $r+1$-th iteration.

In order to compare the invariant $F(G^{\{k\}},t)$ with the $k$-Weisfeiler-Lehman refinement,
 we define a graph invariant $I_{G, k}$ which captures the result of the  $k$-WL coloring and,
  at the same time, it is a combinatorial analogue of $F(G^{k},t)$. For each round $r$,
   we collect all the resulting colors in the multiset
$$
M_{G,\;k}^{r} = \sum_{(i_1...\;i_k)\; \in \mathcal{U}_k} W_{G,\;k}^{r}(i_1...i_k)
$$
Then we define the formal power series
$$
I_{G,\;k}(t)=  \sum_{r=0}^{\infty}t^{r} M_{G,\;k}^{r}
$$

The following technical proposition will be used later.

\begin{prop} \label{Unique permutation}
Let $G$ and $H$ be two graphs with $n$ vertices. Then,
$I_{G,\;k}(t)=I_{H,\;k}(t)$ if and only if there is a permutation $\sigma$
 of the set of $k$-tuples such that
$W_{G,\;k}^{r}(i_1...i_{k})= W_{H,\;k}^{r}(\sigma(i_1...i_{k}) )$ for all $r\geq1$.
 In particular, $$\:tp\:(i_1...i_{k})=\:tp\:(\sigma(i_1...i_{k})).$$
\end{prop}

\textit{Proof:}
The \textit{if} part is immediate. Conversely, assume
$I_{G,\;k}(t)=I_{H,\;k}(t)$.  The coefficient of $t^r$, when
$r=n^k$,  implies
 the existence of  a permutation  $\sigma$ of the set of $k$-tuples  such that
\begin{equation} \label{permutation1}
W_{G,\;k}^{n^k}(i_1...i_{k})= W_{H,\;k}^{n^k}(\sigma(i_1...i_{k}) )
\end{equation}
Whenever Eq. \ref{permutation1} holds for some particular round $r_0$, it holds for all $1 \leq r \leq r_0$. Then,
\begin{equation} \label{permutation2}
W_{G,\;k}^{r}(i_1...i_{k})= W_{H,\;k}^{r}(\sigma(i_1...i_{k}) )
\end{equation}
for all $1 \leq r \leq n^k$. In addition, since the $WL$ refinement stabilizes after the $n^k$ round, we see that Eq.
\ref{permutation2} is true for $r \geq n^k$.
The last assertion is obtained by setting $r=1$ in Eq. \ref{permutation2}. \qed

\section{Graph spectrum is weaker than the $2$-WL refinement}
\label{Graph spectra is weaker than the $2$-WL refinement}

As a warm-up we start by showing that the spectrum of a graph is a
weaker invariant than  the 2-Weisfeiler-Lehman coloring algorithm.
This case displays the essential ingredients of the proof for arbitrary $k$.

\begin{thm}
\label{Case k=1}
 Let $G$ and $H$ be two graphs with  adjacency matrices $A_{G}$ and
$A_{H}$, respectively. If $W_{G,2}^r(i, j) = W_{H,2}^{r}(p,q)$
then $A_{G}^{r}(i,j)=A_{H}^{r}(p,q).$
\end{thm}
\vspace{0.2cm}

\textit{Proof:} We use induction on the number of rounds $r$.
The base case ($r=1$) is trivial.
Assume the statement is valid for $r$, and suppose that
$$W_{G,2}^{r+1}(i, j) = W_{H,2}^{r+1}(p, q).$$
Then, by the definition of the WL coloring,
$$ \sum_{m} (tp_{G}(i, j, m),W_{G,2}^{r}(i, m),W_{G,2}^{r}(m,
j)) =\sum_{m} (tp_{H}(p, q,
m),W_{H,2}^{r+1}(p,m),W_{H,2}^{r+1}(m, q)).$$

This is an equality of multisets. This means that there exists a
permutation $\sigma$ of $ \{ 1,2,...,n \} $ such that

$$
\left\{
\begin{array}{cccc}
     & tp_{G}(i, j, m)&= &tp_{H}(p, q, \sigma(m)),\\
     & W_{G,2}^{r}(i, m)&=&W_{H,2}^{r}(p,\sigma(m)),\\
     & W_{G,2}^{r}(m, j)&=&W_{H,2}^{r}(\sigma(m),q)).\\
\end{array}
\right.
$$

By the induction hypothesis, this implies

$$
\left\{
\begin{array}{ll}
     & A_{G}(i,m)= A_{H}(p,\sigma(m)), \;  \;  A_{G}(m,j)= A_{H}(\sigma(m),q), \\
     & A^{r}_{G}(i,m)=A^{r}_{H}(p,\sigma(m))\\
     & A^{r}_{G}(m, j)=A^{r}_{H}(\sigma(m),q)\\
\end{array}.
\right.
$$

Summing over $m$, we have

 $$ \sum_{m} A_{G}(i,m) A^{r}_{G}(m,j)=
    \sum_{m} A_{H}(p,m) A^{r}_{H}(m,q),$$

that is, $A^{r+1}_{G}(i,j)=A^{r+1}_{H}(p,q)$  \qed

\begin{thm} \label{Main theorem for k=1}
 Let $G$ and $H$ be two graphs. If  $I_{G,\;2}(t)=I_{H,\;2}(t)$,  then $G$ and $H$ are cospectral.
 \end{thm}

\textit{Proof:}
Assume $I_{G,\;2}(t)=I_{H,\;2}(t)$. By Prop. \ref{Unique permutation}, there is a permutation  $\sigma$ of the set of $2$-tuples  such that, for every $2$-tuple $ij$,
$$
W_{G,\;2}^{r}(ij)= W_{H,\;2}^{r}(\sigma(ij) )
$$
for $r \geq 1$. When $r=1$, this is
$$\:tp\:(ij)=\:tp\:(\sigma(ij)).$$
In particular, $\sigma$ sends the diagonal of $W_{G,\;2}^{r}$ to the diagonal of $W_{H,\;2}^{r}$, that is,
$$\sigma (ii)=pp$$ for some element $p$.
Then, collecting all the colors in the diagonal, we have
$$
\sum_{i} W_{G,\;2}^{r}(ii)=\sum_{i} W_{H,\;2}^{r}(\sigma(i)\sigma(i))
$$
By  Theorem \ref{Case k=1}, this implies
$$
\sum_{i} A^{r}_{G}(i,i)=\sum_{i} A^{r}_{H}(\sigma(i),\sigma(i))
$$
that is, $Tr A^{r}_{G} = Tr A^{r}_{H}$ for $r \geq 1$. Then, $F(G,t)=F(H,t)$ and this means that $G$ and $H$ are cospectral.
  \qed

\section{Spectra of $k$-th powers}
\label{Spectrum of the $k$-th power vs.  the $2k$-WL coloring}

For each round $r$, we think of the  $2k$-WL coloring as a matrix of colors: the rows and columns are indexed by $k$-tuples, with the color $W_{G,\;k}^{r}(i_1...i_k j_1...j_k)$ in the entry  $(i_1...i_k , j_1...j_k)$.

\begin{thm} \label{core}
 Let $G^{k}$ and $H^{k}$ be the $k$-th powers of two graphs $G$ and $H$ respectively. Let $A_{G^{k}}^{r}$ and $A_{H^{k}}^{r}$ be the $r$-th powers of their adjacency matrices. If
 $$W_{G,2k}^{r}(i_1\dots i_k \ j_1\dots j_k) = W_{H,2k}^{r}(p_1\dots p_k \ q_1\dots q_k),$$
 then
$$A_{G^{k}}^{r}(i_1 \dots i_k, j_1 \dots j_k)=A_{H^{k}}^{r}(p_1 \dots p_k, q_1 \dots q_k).$$
\end{thm}

\textit{Proof:} The proof goes along the lines of Theorem \ref{Case k=1}. Let
$r=1$. Suppose that $$A_{G^{k}}(i_1 \dots i_k, j_1 \dots j_k)=1.$$
Then $i_{l}=j_{l}$ for all $l$ except for a unique value $l_{0}$,
for which $A_{G}(i_{l_0},j_{l_0})=1$. By hypothesis,
$$W_{G,2k}^{1}(i_1 \dots i_k \ j_1\dots j_k) = W_{H,2k}^{1}(p_1 \dots p_k \ q_1 \dots q_k),$$
that is,
$$tp(i_1 \dots i_k \ j_1\dots j_k) = tp(p_1 \dots p_k \ q_1 \dots q_k).$$
By the definition of type, this implies that
$p_{l}=q_{l}$ for $l \neq l_{0}$ and $A_{H}(p_{l_0},q_{l_0})=1$.
 Then
$A_{H^{k}}(p_1\dots p_k, q_1\dots q_k)=1$. The argument can be
reversed, proving that $$A_{G^{k}}(i_1\dots i_k, j_1 \dots
j_k)=A_{H^{k}}(p_1\dots p_k, q_1\dots q_k).$$ This prove the case
$r=1$. Now assume the statement is valid
for $r$, and suppose that
$$W_{G,2k}^{r+1}(i_1 \dots i_k \  j_1 \dots j_k) = W_{H,2k}^{r+1}(p_1 \dots p_k \ q_1 \dots q_k).$$
By the definition of the WL coloring,

$$ \sum_{m \in V} \left( tp_{G}(i_1 \dots i_k \ j_1 \dots j_k \  m),
S_{G,2k}^{r}(i_1\dots i_k \  j_1\dots j_{k}\ m)\right)=$$
$$=\sum_{m \in V} \left( tp_{H}(p_1 \dots p_k \ q_1 \dots
q_k \  m), S_{H,2k}^{r}(p_1\dots p_k \  q_1\dots q_{k}\
m)\right).  $$
Therefore  there exists a permutation $\sigma$ of $ \{ 1,2,...,n
\} $ such that

$$
\left\{
\begin{array}{cccc}
     & tp_{G}(i_1\dots i_k \ j_1\dots j_k \  m)&=& tp_{H}(p_1\dots p_k \ q_1\dots q_k \  \sigma(m)),\\
     & W_{G,2k}^{r}(i_1\dots i_k \ j_1\dots j_{k-1} \  m)&=&W_{H,2k}^{r}(p_1\dots p_k \ q_1 \dots q_{k-1}, \sigma(m)),\\
     & &\dots & \\
     & W_{G,2k}^{r}(m \ i_2\dots i_k \ j_1\dots j_{k})&=&W_{H,2k}^{r}(\sigma(m) \ p_2 \dots p_k \ q_1 \dots q_{k}).\\
\end{array}
\right.
$$
The induction hypothesis implies
$$\left\{
\begin{array}{ll}
& A_G(i_t,m)=A_G(p_t,\sigma(m)) \quad \quad \mbox{for } t=1,\dots, k.\\
& A^{r}_{G^k}(i_1\dots m \dots i_k , j_1\dots
j_k)=A^{r}_{H^k}(p_1\dots \sigma(m)\dots  p_k , q_1\dots q_k ).
\end{array}
\right.
$$

Our goal is to show that

$$
A^{r+1}_{G^k}(i_1...i_k, j_1...j_k)=A^{r+1}_{H^k}(p_1...p_k,
q_1...q_k).
$$

We have
\begin{equation}
A^{r+1}_{G^k}(i_1...i_k, j_1...j_k) =
\sum_{s_1...s_k}A_{G^k}(i_1...i_k, s_1...s_k) A^{r}_{G^k}(s_1...s_k, j_1...j_k)
\end{equation}
Observe that  $A_{G^k}(i_1...i_k, s_1...s_k)=0$ unless there
exists an index $t$ such that $A_G(i_t,s_t)=1$ and $i_l=s_l$ for
all $l \neq t$. Hence
$$
A^{r+1}_{G^{k}}(i_1\dots i_k, j_1\dots j_k) = \sum_{m \in V}
\sum_{t=1}^k A_{G}(i_t, m) A^{r}_{G^{k}}(i_1 \dots m \dots i_k, j_1
\dots j_k)
$$
$$
=\sum_{m \in V} \sum_{t=1}^k A_{H}(p_t, \sigma(m))
A^{r}_{H^k}(p_1 \dots \sigma(m) \dots p_k, q_1
\dots q_k)
= A^{r+1}_{H^k}(p_1\dots p_k, q_1\dots q_k)
$$
   \qed

\begin{thm} \label{Main theorem for powers}
 Let $G$ and $H$ be two graphs. If  $I_{G,\;2k}(t)=I_{H,\;2k}(t)$, then $$F(G^k,t)=F(H^k,t).$$ In other words, if the $2k$-th WL refinement cannot distinguish $G$ from $H$, then their $k$-th powers are cospectral.
 \end{thm}

\textit{Proof:}
Assume $I_{G,\;2k}(t)=I_{H,\;2k}(t)$. By Prop. \ref{Unique permutation}, there is a permutation  $\sigma$ of the set of $2k$-tuples  such that, for every $2k$-tuple $i_1...i_{k}j_1...j_{k}$,
$$
W_{G,\;2k}^{r}(i_1...i_{k}j_1...j_{k})= W_{H,\;2k}^{r}(\sigma(i_1...i_{k}j_1...j_{k}) )
$$
for $r \geq 1$. When $r=1$, this is
$$\:tp\:(i_1...i_{k}j_1...j_{k})=\:tp\:(\sigma(i_1...i_{k}j_1...j_{k}).$$
In particular, $\sigma$ sends the diagonal of $W_{G,\;2k}^{r}$ to the diagonal of $W_{H,\;2k}^{r}$, that is,
$$\sigma (i_1...i_{k}i_1...i_{k})=p_1...p_{k}p_1...p_{k}$$ for some $k$-tuple $p_1...p_{k}$.
Then, collecting all the colors in the diagonal, we have
$$
\sum_{i_1...i_{k}} W_{G,\;2k}^{r}(i_1...i_{k}i_1...i_{k})=\sum_{i_1...i_{k}} W_{H,\;2k}^{r}(\sigma(i_1...i_{k})\sigma(i_1...i_{k}))
$$
By Theorem \ref{core}, this implies
$$
\sum_{i_1...i_{k}} A^{r}_{G^k}(i_1...i_{k},i_1...i_{k})=\sum_{i_1...i_{k}} A^{r}_{H^k}(\sigma(i_1...i_{k}),\sigma(i_1...i_{k}))
$$
that is, $Tr A^{r}_{G^k} = Tr A^{r}_{H^k}$ for $r \geq 1$. Then, $F(G^k,t)=F(H^k,t)$.
  \qed

Our goal is to prove the analogue of Theorem \ref{Main theorem for powers} for the $k$-th symmetric powers. As an intermediate step, we prove analogues of Theorem \ref{core} and Theorem \ref{Main theorem for powers} for the restricted $k$-th powers.

\begin{thm} \label{restricted version}
 Let $G^{(k)}$ and $H^{(k)}$ be the $k$-th restricted powers of two graphs $G$ and $H$. Let $A_{G^{(k)}}^{r}$ and $A_{H^{(k)}}^{r}$ be the $r$-th powers of their adjacency matrices.
 Assume that $i_1\dots i_k$, $j_1\dots j_k$, $p_1\dots p_k$, and $q_1\dots q_k$ are $k$-tuples in $\mathcal{D}_k$.
 If
 $$W_{G,2k}^{r}(i_1\dots i_k \ j_1\dots j_k) = W_{H,2k}^{r}(p_1\dots p_k \ q_1\dots q_k),$$
 then
$$A_{G^{(k)}}^{r}(i_1 \dots i_k, j_1 \dots j_k)=A_{H^{(k)}}^{r}(p_1 \dots p_k, q_1 \dots q_k).$$
\end{thm}

\textit{Proof:}
The proof mimics that of Theorem \ref{core}. The case $r=1$ is unaltered, so we assume the proposition is valid for $r$
and we suppose that
$$W_{G,2k}^{r+1}(i_1 \dots i_k \  j_1 \dots j_k) = W_{H,2k}^{r+1}(p_1 \dots p_k \ q_1 \dots q_k).$$
This means that there is a permutation
$\sigma$ of $ \{ 1,2,...,n
\} $ such that
$$
\left\{
\begin{array}{cccc}
     & tp_{G}(i_1\dots i_k \ j_1\dots j_k \  m)&=& tp_{H}(p_1\dots p_k \ q_1\dots q_k \  \sigma(m)),\\
     & W_{G,2k}^{r}(i_1\dots i_k \ j_1\dots j_{k-1} \  m)&=&W_{H,2k}^{r}(p_1\dots p_k \ q_1 \dots q_{k-1}, \sigma(m)),\\
     & &\dots & \\
     & W_{G,2k}^{r}(m \ i_2\dots i_k \ j_1\dots j_{k})&=&W_{H,2k}^{r}(\sigma(m) \ p_2 \dots p_k \ q_1 \dots q_{k}).\\
\end{array}
\right.
$$
From the first of these equations, we observe that $m = i_t$ implies $\sigma(m) = p_t$. Therefore, the $k$-tuple
 $(i_1 ...i_{l-1} \; m \; i_{l+1} ... i_k)$ is in $\mathcal{D}_k$ if and only if  $$(p_1 ...p_{l-1} \; \sigma(m) \; p_{l+1} ... p_k)$$ is  in $\mathcal{D}_k $.

This observation shows that, if we assume $m \neq i_t$ for $t=1,...,k$, we are allowed to apply the induction hypothesis to obtain
$$\left\{
\begin{array}{ll}
& A_G(i_t,m)=A_G(p_t,\sigma(m)) \quad \quad \mbox{for } t=1,\dots, k.\\
& A^{r}_{G^k}(i_1\dots m \dots i_k , j_1\dots
j_k)=A^{r}_{H^k}(p_1\dots \sigma(m)\dots  p_k , q_1\dots q_k ).
\end{array}
\right.
$$
Then
\begin{equation}
A^{r+1}_{G^{(k)}}(i_1...i_k, j_1...j_k) =
\sum_{(s_1...s_k) \in \mathcal{D}_k}A_{G^{(k)}}(i_1...i_k, s_1...s_k) A^{r}_{G^{(k)}}(s_1...s_k, j_1...j_k)
\end{equation}
$$
= \sum_{m \notin \{i_1,...,i_k\}} \sum_{t=1}^k A_{G}(i_t, m) A^{r}_{G^{(k)}}(i_1 \dots m \dots i_k, j_1
\dots j_k)
$$
$$
= \sum_{\sigma (m) \notin \{p_1,...,p_k\}} \sum_{t=1}^k A_{H}(p_t, \sigma(m))
A^{r}_{H^{(k)}}(p_1 \dots \sigma(m) \dots p_k, q_1
\dots q_k)
$$
$$
= A^{r+1}_{G^{(k)}}(p_1\dots p_k, q_1\dots q_k)
\qed $$

\begin{thm} \label{Main theorem for restricted powers}
 If the $2k$-th WL refinement fails to distinguish  $G$ from $H$, then their restricted $k$-th powers are cospectral.
 \end{thm}

\textit{Proof:}
The proof is analogous to that of Theorem \ref{Main theorem for powers}. Assume $I_{G,\;2k}(t)=I_{H,\;2k}(t)$.
Let $\sigma$ be the permutation of the set of $2k$-tuples given by Proposition \ref{Unique permutation}. Since $\sigma$ preserves the type of the $2k$-tuples, if
$i_1...i_{k}$ is in $\mathcal{D}_{k}$, then
$$\sigma (i_1...i_{k}i_1...i_{k})=p_1...p_{k}p_1...p_{k}$$ for some $k$-tuple $p_1...p_{k}$ $\in \mathcal{D}_{k}$.
Then,
$$
\sum_{(i_1...i_{k}) \in \mathcal{D}_{k}} W_{G,\;2k}^{r}(i_1...i_{k}i_1...i_{k})=\sum_{(i_1...i_{k}) \in \mathcal{D}_{k}} W_{H,\;2k}^{r}(\sigma(i_1...i_{k})\sigma(i_1...i_{k}))
$$
By Theorem \ref{restricted version}, this implies
$$
\sum_{(i_1...i_{k}) \in \mathcal{D}_{k}} A^{r}_{G^{(k)}}(i_1...i_{k},i_1...i_{k})=\sum_{(i_1...i_{k}) \in \mathcal{D}_{k}} A^{r}_{H^{(k)}}(\sigma(i_1...i_{k}),\sigma(i_1...i_{k}))
$$
that is, $Tr A^{r}_{G^{(k)}} = Tr A^{r}_{H^{(k)}}$ for $r \geq 1$. Then, $F(G^{(k)},t)=F(H^{(k)},t)$.
  \qed

\section{Proof of Theorem 1}
\label{Proof of Theorem 1}

We can restate Theorem \ref{Main theorem} as follows:

 \begin{thm}
 Let $G$ and $H$ be two graphs. If  $I_{G,\;2k}(t)=I_{H,\;2k}(t)$, then $$F(G^{\{k\}},t)=F(H^{\{k\}},t).$$
 \end{thm}

\textit{Proof:}
Assume $I_{G,\;2k}(t)=I_{H,\;2k}(t)$. Again, by Prop. \ref{Unique permutation}, there is a permutation  $\sigma$ of the set of $2k$-tuples  such that
\begin{equation} \label{hypothesis}
W_{G,\;2k}^{r}(i_1...i_{2k})= W_{H,\;2k}^{r}(\sigma(i_1...i_{2k}) )
\end{equation}
for all $r \geq 1$. Since
$$\:tp\:(i_1...i_{2k})=\:tp\:(\sigma(i_1...i_{2k})),$$
we can restrict $\sigma$ in the following way.
 If $\theta$ is a permutation in  $S_k$, we denote by $\theta(i_1...i_{k})$ the $k$-tuple $(i_{\theta(1)}...i_{\theta(k)})$.
 Let us write the   $2k$-tuples as pairs of $k$-tuples: $(i_1...i_{k},j_1...j_{k})$.
 Observe that if a $2k$-tuple is of the form  $$(i_1...i_{k},\theta(i_1...i_{k}))$$,where   $(i_1...i_{k}) \in \mathcal{D}_k$ and  $\theta \in S_k$, 
  then (due to the type-conservation) $\sigma$ sends it to a $2k$-tuple of the form $(j_1...j_{k},\theta(j_1...j_{k}))$, for some $(j_1...j_{k}) \in \mathcal{D}_k$. Thus, there is a permutation $\omega$ of the set $\mathcal{D}_k$ such that for every $(i_1...i_{k}) \in \mathcal{D}_k$

\begin{equation}
W_{G,\;2k}^{r}(i_1...i_{k},\theta(i_1...i_{k}))= W_{H,\;2k}^{r}(\omega(i_1...i_{k}),\theta (\omega (i_1...i_{k})) )
\end{equation}
By Theorem \ref{restricted version}, it follows that
\begin{equation}
A_{G^{(k)}}^{r}(i_1...i_{k},\theta(i_1...i_{k}))= A_{H^{(k)}}^{r}(\omega(i_1...i_{k}),\theta (\omega (i_1...i_{k})) ).
\end{equation}
In particular,
\begin{equation}
\sum_{(i_1...i_{k}) \in \mathcal{D}_k}\sum_{\theta \in S_k} A_{G^{(k)}}^{r}(i_1...i_{k},\theta(i_1...i_{k}))= \sum_{(i_1...i_{k}) \in \mathcal{D}_k}\sum_{\theta \in S_k} A_{H^{(k)}}^{r}(\omega(i_1...i_{k}),\theta (\omega (i_1...i_{k})) ).
\end{equation}
Since $\omega$ is a bijection, we can drop it from this last equation, and  we have
\begin{equation}
\sum_{(i_1...i_{k}) \in \mathcal{D}_k}\sum_{\theta \in S_k} A_{G^{(k)}}^{r}(i_1...i_{k},\theta(i_1...i_{k}))= \sum_{(i_1...i_{k}) \in \mathcal{D}_k}\sum_{\theta \in S_k} A_{H^{(k)}}^{r}(i_1...i_{k},\theta  (i_1...i_{k}) )
\end{equation}
Let $M^k$ be the matrix of Prop. \ref{Quotient}. This  last equation can be written as
$$
Tr (A_{G^{(k)}}^{r}M_k) =  Tr (A_{H^{(k)}}^{r}M_k)
$$
By Prop. \ref{Quotient}, this is equivalent to
\begin{equation}
Tr (A_{G^{\{k\}}}^{r}) = Tr (A_{H^{\{k\}}}^{r})
\end{equation}
Since this is true for all $r$, then
$F(G^{\{k\}},t)=F(H^{\{k\}},t).$




\end{document}